\documentclass[11pt]{article}
\title{A characterization property on fields equivalent to algebraicity on Banach spaces}
\author{Xavier Le Breton}
\date{}
\usepackage{amssymb}
\usepackage{epsf}
\usepackage{fullpage}
\usepackage{theorem}
\newtheorem{theorem}{Theorem}

\newtheorem{corollary}[theorem]{Corollary}
\newtheorem{proposition}[theorem]{Proposition}

\newtheorem{definition}[theorem]{Definition}

\def\labibliographie#1{\section*{References}\list
 {[\arabic{enumi}]}{\settowidth\labelwidth{[#1]}\leftmargin\labelwidth
 \advance\leftmargin\labelsep
 \usecounter{enumi}}
 \def\newblock{\hskip .11em plus .33em minus -.07em}
 \sloppy
 \sfcode`\.=1000\relax}

\begin{document}
\maketitle
\begin{center}
Laboratoire de Recherche en Informatique \\
B\^atiment $490$ Universit\'e Paris-Sud \\
$91405$ Orsay CEDEX FRANCE \\
lebreton@lri.fr
\end{center}

In $1980$, Christol, Kamae, Mend\`es France and Rauzy stated in \cite{ckmr} an important theorem in automata theory. This theorem links sequences recognized by automata and algebraic formal power series. In $1994$, Bruy\`ere, Hansel, Michaux and Villemaire extended this theorem with a logical link in \cite{bhmv}. With theses two articles, we can translate the property for a formal power series to be algebraic in combinatorics terms or logical terms. Our general purpose is to extend these theorems to algebraic dependences between formal power series. We want to be able to translate in combinatorial terms the fact, for two formal power series, to be algebraically dependent. Our first approach, see \cite{lb}, was combinatorial, and we proved linear independences between some formal power series. The second idea is to use logic (remember that for the case of algebraicity theses points of view are equivalent) and hope it could be translated in combinatorial terms further. That is why we were interested in the work of Tyszka (even if it does not speak about formal power series). Indeed, Tyszka introduce a logical property which is equivalent to algebraicity in $\mathbb{R}$ and in the $p$-adic field $\mathbb{Q}_p$. The goal of this article is to study this property and describe fields for which it is equivalent to algebraicity. We will see that the formal power series field is one of them and why finding a good equivalence for algebraic dependence is not easy. Actually, the source of the problem is that we work on a field with positive characteristic so we suggest a property quite different from algebraic dependence but (we hope) more likely equivalent to some combinatorial characterization.

\section{General backgrounds}
In this section, we present the general background needed to
understand other sections. We begin with a presentation of the
work of Tyszka in \cite{ty} which is the base of this article.
Next we define the notion of automatic sequences and talk about
the theorem of Christol. At last, we recall some logical material.

\subsection{Tyzska's article}
In \cite{ty}, Tyzska studies an interesting property which
characterizes elements in a field. In order to cleary explain how
it works we will define a family of functions we will call
pseudo-morphisms. A pseudo-morphism tries to preserve ``as much as
possible'' the algebraic structure of its domain.

\medskip

\begin{definition}
Let $(K,0,1,+,\times)$ be a field and A be a subset of K. We will
say that $\phi$ : $A\rightarrow K$ is a pseudo-morphism if and
only if it satisfies:
\begin{itemize}
\item $0\!\in\! A\Rightarrow \phi (0)=0$
\item $1\!\in\! A\Rightarrow \phi (1)=1$
\item $a,b,(a+b)\!\in\! A \Rightarrow \phi (a) + \phi(b) = \phi(a+b)$
\item $a,b,(a\times b)\!\in\! A \Rightarrow \phi (a) \times \phi(b) = \phi(a\times b)$
\end{itemize}
\end{definition}

\medskip

Hence, a pseudo-morphism preserves all valid additions and
multiplications within its domain. Notice that a pseudo-morphism
with a domain which is a field is a field morphism. The smaller
the domain is, the easier it is to find a pseudo-morphism. For
example, any function from $\{2,5\}$ to $\mathbb{R}$ is a
pseudo-morphism. However, a pseudo-morphism $\phi$ from
$\{2,3,5\}$ to $\mathbb{R}$ must satisfy
$\phi(2)+\phi(3)=\phi(5)$.

\medskip

\begin{proposition}
Let $\phi$ : $A\rightarrow K$ be a pseudo-morphism and $A'\subset
A$, then $\phi_{\vert_{A'}}$ : $A'\rightarrow K$ is also a pseudo
morphism.
\end{proposition}

\medskip

We can now present the property studied in \cite{ty}, it deals
with pseudo-morphims having a finite domain.

\medskip

\begin{definition}
Let $(K,0,1,+,\times)$ be a field and $x$ an element in $K$, $x$
is Tyzska characterizable $(TC)$ if and only if there exists a
finite set $A$ included in $K$ and containing $x$ so that:
$$\forall \phi \textrm{ : } A\rightarrow K \textrm{ , } \phi \textrm{ is a pseudo-morphism }\Rightarrow \phi(x)=x.$$
\end{definition}

\medskip

An element is characterizable if it gets a unique behavior (in
respect to addition and mutiplication) within a finite set. We
give a proposition proved in \cite{ty}.

\medskip

\begin{proposition}
The set (noted $\widetilde{K}$) of all characterizable elements in
$K$ is always a subfield of $K$.
\end{proposition}

\medskip

In \cite{ty}, Tyszka studied this characterization in classical
fields of characteristic $0$. He proved the following results:

\medskip

\begin{itemize}
\item $\widetilde{\mathbb{R}}=\{x\!\in\! \mathbb{R}\textrm{ :
}x\textrm{ is algebraic over }\mathbb{Q}\}$ \item For all $p$
prime number, $\widetilde{\mathbb{Q}_p}=\{x\!\in\!
\mathbb{Q}_p\textrm{ : }x\textrm{ is algebraic over }\mathbb{Q}\}$
\item $\widetilde{\mathbb{C}}=\mathbb{Q}$
\end{itemize}

\medskip

The proof of the first two results uses a theorem of model theory.
In this article, we present a general theorem for fields which are
also Banach algebras with a proof using analytic methods. The
third result above is astonishing: we cannot prove a theorem
saying that if a field is also a $\mathbb{Q}$-Banach
($\mathbb{R}$, $\mathbb{Q}_p$ and $\mathbb{C}$ are this kind of
fields) then the Tyszka property characterizes all
$\mathbb{Q}$-algebraic elements. We will have to define another
property to characterize all algebraic elements in $\mathbb{C}$
(not only the rational ones).

\medskip

In fact, there is another reason to change the Tyszka property. It
is not hard to prove that for any field $K$, any element in the
prime field of $K$ is characterizable. The problem in $\mathbb{C}$
is that only elements in its prime field are characterizable. Our
purpose is to study the Tyszka property on a typical field of
positive characteristic: $\mathbb{F}_p((X))$. However, Tyszka also
proved a general theorem for this kind of fields:

\medskip

\begin{proposition}
Let $K$ be a field of positive characteristic then $\widetilde{K}$
is the prime field of $K$.
\end{proposition}

\medskip

Moreover, the case of the field $\mathbb{C}$ is not singular, as
Tyszka proved in \cite{ty2}.

\medskip

\begin{proposition}
Let $K$ be a field, if there exists a subfield of $K$ which is
algebraically closed then $\widetilde{K}$ is the prime field of
$K$
\end{proposition}

\medskip
Now, because $\mathbb{C}$ is an algebraically closed field, the
fact that $\widetilde{\mathbb{C}}=\mathbb{Q}$ is not so
surprising. We will explain later how to find a Tyszka-like
property that can also work on algebraically closed fields.

\subsection{Automatic sequences}
Now, we give a short presentation of automatic sequences. The
reader may consult \cite{alsh} and \cite{ckmr}. An automatic
sequence is a sequence generated by a deterministic finite
automaton with output (DFAO).

\medskip

\begin{definition}
A DFAO is a $6$-uplet $M=(Q,\Sigma,\delta,q_0,\Delta,\tau)$ where
\\ $Q$ is the finite set of states,
\\ $\Sigma$ is a finite input alphabet,
\\ $\delta$ : $Q\times\Sigma\rightarrow Q$ is the transition function,
\\ $q_0\!\in\! Q$ is the initial state,
\\ $\Delta$ is the output alphabet, and
\\ $\tau$ : $Q\rightarrow\Delta$ is the transition function.
\end{definition}

\medskip

Each DFAO $M$ defines a function from $\Sigma^\ast$ to $\Delta$
denoted $f_M$. Let $w=w_1w_2\cdots w_k\!\in\!\Sigma^\ast$ where
$k=\vert w\vert$. If $k=0$ then $w=\varepsilon$ and we define
$f_M(w)=\tau(q_0)$. Suppose$\omega\not=\epsilon$, the automaton
first reads the letter $w_1$ being in its initial state $q_0$: it
moves from $q_0$ to $q_1=\delta(q_0,w_1)$. Then the automaton
continues by reading $w_2$ while being in the state $q_1$: it
moves from $q_1$ to $q_2=\delta(q_1,w_2)$. The automaton continues
to read all letters until the last one and stops in the state
$q_k=\delta(q_{k-1},w_k)$. We define $f_M$ by applying the output
function to this state.

\medskip

$$f_M(w)=\tau(q_k)=\tau(\delta(\cdots \delta(\delta(q_0,w_1),w_2)\cdots,w_k))$$

\medskip

Let $p$ be a prime number, we are now able to define $p$-automatic
sequences. We consider all automata
$M=(Q,\Sigma,\delta,q_0,\Delta,\tau)$ where
$\Sigma=\Delta=\mathbb{F}_p$. This kind of automata defines
functions from $\mathbb{F}_p^\ast$ to $\mathbb{F}_p$. But the set
$\mathbb{F}_p^\ast$ is in bijection with the set of natural
numbers $\mathbb{N}$. In fact, each $w\!\in\!\mathbb{F}_p^\ast$ is
the $p$-base expansion of a natural number, where $\varepsilon$ is
the $p$-base expansion of $0$. By this way the function $f_M$
becomes a sequence over $\mathbb{F}_p$.

\medskip

\begin{definition}
Let $p$ be a prime number. The sequence
$u=(u_n)_{n\!\in\!\mathbb{N}}$ over $\mathbb{F}_p$ is called a
$p$-automatic sequence if and only if there is an automaton
$M=(Q,\mathbb{F}_p,\delta,q_0,\mathbb{F}_p,\tau)$ generating $u$.
\end{definition}

\subsection{Theorem of Christol and Cartier's operators}
The set of all sequences over $\mathbb{F}_p$ is naturally in
bijection with the set of formal power series with coefficients in
$\mathbb{F}_p$. This bijection maps a sequence
$(u_n)_{n\!\in\!\mathbb{N}}$ to the formal power series
$\Sigma_{n\!\in\!\mathbb{N}}u_nX^n$. The set $\mathbb{F}_p[[X]]$
is actually a ring. An element of this ring is called algebraic if
it is a zero of a non-trivial polynomial with coefficients in
$\mathbb{F}_p[X]$. Christol's theorem says that algebraicity and
$p$-automaticity are equivalent from this point of view.

\medskip

\begin{theorem}
Let $u=(u_n)_{n\!\in\!\mathbb{N}}$ and
$F(X)=\sum_{n\!\in\!\mathbb{N}}u_nX^n\!\in\!\mathbb{F}_p[[X]]$ its
associated formal power series. Then $u$ is $p$-automatic if and
only if $F(X)$ is algebraic.
\end{theorem}

\medskip

The proof is in \cite{ckmr}. We just give the main point of the
proof which permits to introduce Cartier's operators.

\medskip

\begin{definition}
We define $\Lambda_0,\Lambda_1,\cdots,\Lambda_{p-1}$ functions
from $\mathbb{F}_p[[X]]$ to itself as follows. Let
$F(X)=\Sigma_{n\!\in\!\mathbb{N}}u_nX^n$ and $i\!\in\![0,p-1]$.
$$\Lambda_i(F(X))=\Sigma_{n\!\in\!\mathbb{N}}u_{pn+i}X^n$$
We call $\Lambda_0,\Lambda_1,\cdots,\Lambda_{p-1}$ the Cartier operators.
\end{definition}

\medskip

The first thing to notice on these operators is that for every
formal power series $F(X)\!\in\!\mathbb{F}_p[[X]]$, the following
equation stands:
$$F(X)=\Sigma_{i=0}^{p-1}X^i(\Lambda_i(F(X)))^p$$
Furthermore, for every
$Z=(Z_0(X),\cdots,Z_{p-1}(X))\!\in\!(\mathbb{F}_p[[X]])^p$ we
have:
$$F(X)=\Sigma_{i=0}^{p-1}X^i(Z_i(X))^p\Rightarrow Z=(\Lambda_0(F(X)),\cdots,\Lambda_{p-1}(F(X)))$$
The Cartier operators can be used to define the $p$-kernel of a formal power series.

\medskip

\begin{definition}
Let $F(X)\!\in\!\mathbb{F}_p[[X]]$. We define the $p$-kernel
$K_p(F(X))$ as the smallest set (with respect to inclusion)
containing $F(X)$ and preserved by Cartier's operators.
\end{definition}

\medskip

The idea of the proof of theorem of Christol is to prove the
following two equivalences:
$$ u\textrm{ is $p$-automatic }\Leftrightarrow K_p(F(X))\textrm{ is a finite set.}$$
$$ F(X)\textrm{ is algebraic }\Leftrightarrow K_p(F(X))\textrm{ is a finite set.}$$

\subsection{Some logical material}
In this paragraph we recall definitions of a logical structure and
a first-order formula.

\medskip

\begin{definition}
We call $S=\{D,(R_i)_{i\!\in\! I},(f_j)_{j\!\in\!
J},(c_k)_{k\!\in\! K} \}$ a logical structure when $D$ is a set,
$(R_i)_{i\!\in\! I}$ is a family of relations on $D$ (that is a
subset of $D^{n_i}$), $(f_j)_{j\!\in\! J}$ is a family of
functions on $D$ (from $D^{n_j}$ to $D$) and $(c_k)_{k\!\in\! K}$
is a family of elements in $D$ (called constants).
\end{definition}

\medskip

For example $(\mathbb{R},+,\times ,0,1)$ is a logical structure.
In order to construct the set of first-order formulas we have to
define the set of terms. We first need to fix an infinite subset
$V={x,y,z,t,\cdots}$ called the set of variables.

\medskip

\begin{definition}
We define the set of terms by induction using three rules: \\
$1$ Each constant is a term. \\
$2$ Each variable is a term. \\
$3$ If $j\!\in\! J$ and $t_1,\cdots,t_{n_j}$ are terms, then $f_j(t_1,\cdots,t_{n_j})$ is a term.
\end{definition}

\medskip

\begin{definition}
We define the set of first-order formulas by induction using four rules: \\
$1$ If $t_1$ and $t_2$ are terms, then $(t_1=t_2)$ is a formula. \\
$2$ If $i\!\in\! I$ and $t_1,\cdots,t_{n_i}$ are terms, then $(R_i(t_1,\cdots,t_{n_i}))$ is a formula. \\
$3$ If $\phi$ and $\psi$ are formulas, then $(\phi\lor\psi)$,$(\phi\land\psi)$,$(\lnot\phi)$,$(\phi\Rightarrow\psi)$ and $(\phi\Leftrightarrow\psi)$ are formulas. \\
$4$ If $\phi$ is a formula and $x\!\in\! V$, then $(\forall x\phi)$ and $(\exists x\phi)$ are formulas.
\end{definition}

\medskip

A variable $x$ is called \textit{free} in a formula $\phi$ if and
only if it occurs at least one time in $\phi$ and is not under the
scope of a quantifier. A \textit{sentence} is a formula with no
free variable (that is each variable that occurs at least one time
in the formula is under the scope of a quantifier). The set of all
sentences does not depend on the set $D$ in the structure. However
the logical value of a sentence (is it true or false?) is closely
related to the set $D$. For example, let $\phi$ be the following
sentence: $\phi=\exists x,x\times x=1+1$. Then $\phi$ is false in
the logical structure $\{\mathbb{Q},+,\times,0,1\}$ but is true in
$\{\mathbb{R},+,\times,0,1\}$. We say that the set $\mathbb{R}$
valids $\phi$ and we note $\mathbb{R}\vdash\phi$. The set of all
valid sentences of a set $D$ is called the theory of $D$ and noted
${Th}\{D\}$ or ${Th}\{D,(R_i)_{i\!\in\! I},(f_j)_{j\!\in\!
J},(c_k)_{k\!\in\! K}\}$. Tyszka used in his article a
model-theory theorem: the fact that
${Th}\{\mathbb{R}\}={Th}\{\mathbb{R}_{alg}\}$ where
$\mathbb{R}_{alg}$ is the set of all algebraic reals. In our
article we give an analytic proof of the same result.

\section{Fields of characteristic zero}
This section presents and studies two properties ($(FTC)$ and
$(TC)$) in fields of characteristic zero. We prove an equivalence
theorem in Banach algebras and show counterexamples in non-Banach
algebras.

\subsection{Definitions and links with algebraicity}
In this section, we first define two general properties ($(TC)$
and $(FTC)$) looking like the one studied by Tyszka. We then prove
some fact relating algebraicity and theses two properties.

\subsubsection{Definitions}
We want to give a general definition of the property of Tyszka on
logical structures. As in Section $1.1$, we first define
pseudo-morphisms on structures.

\medskip

\begin{definition}
Let $S=\{D,(R_i)_{i\!\in\! I},(f_j)_{j\!\in\! J},(c_k)_{k\!\in\!
K} \}$ be a structure. Let $A\subset D$ and $\phi : A\rightarrow
D$. We say that $\phi$ is a pseudo-morphism on S if and only if:
$$\forall k\!\in\! K, \  c_k\!\in\! A\Rightarrow \phi(c_k)=c_k$$
$$\forall j\!\in\! J, \  n\!\in\!\mathbb{N} \textrm{ so that $f_j$ is an $n$-ary function,}$$
$$\forall (a_1,\cdots ,a_n)\!\in\! A^n, \  f(a_1,\cdots ,a_n)\!\in\! A\Rightarrow \phi(f(a_1,\cdots ,a_n))=f(\phi(a_1),\cdots ,\phi(a_n))$$
$$\forall i\!\in\! I, \  n\!\in\!\mathbb{N} \textrm{ so that $R_i$ is an $n$-ary relation,}$$
$$\forall (a_1,\cdots ,a_n)\!\in\! A^n, \  R(a_1,\cdots ,a_n)\Rightarrow R(\phi(a_1),\cdots ,\phi(a_n))$$
\end{definition}

\medskip

Now, the definition of the generalized Tyszka property is nearly
the same as in Section $1.1$.

\medskip

\begin{definition}
Let $S=\{D,(R_i)_{i\!\in\! I},(f_j)_{j\!\in\! J},(c_k)_{k\!\in\!
K} \}$ be a structure. For every $x\!\in\! D$, we say that $x$ is
Tyszka characterizable $(TC)$ if and only if there exists a finite
set $A\subset D$ such that every $S$-pseudo-morphism $\phi :
A\rightarrow D$ satisfies $\phi(x)=x$.
\end{definition}

\medskip

As we saw in Section $1.1$, the $(TC)$ property cannot be proved
for any algebraic elements of a Banach algebra. We need to define
another property in order to prove a general theorem on Banach
algebras.

\medskip

\begin{definition}
Let $S=\{D,(R_i)_{i\!\in\! I},(f_j)_{j\!\in\! J},(c_k)_{k\!\in\!
K} \}$ be a structure. For every $x\!\in\! D$, we say that $x$ is
Finitely Tyszka characterizable $(FTC)$ if and only if there exist
two finite sets $A\subset D$ and $B\subset D$ so that every
$S$-pseudo-morphism $\phi : A\rightarrow D$ satisfies
$\phi(x)\!\in\! B$.
\end{definition}

\medskip

The idea of the definition starts with the following remark: in
$\mathbb{C}$ the main problem is that we cannot distinguish (using
$(TC)$) all conjugates of an element $x$. For each element in a
field, the set of its conjugates is finite. That is why we replace
$\phi(x)=x$ by $\phi(x)\!\in\! B$ with $B$ a finite subset of $D$.
With this new property we will prove, under reasonable conditions,
the equivalence of being algebraic and being $(FTC)$.

\subsubsection{Algebraicity, $(TC)$ and $(FTC)$}
Given a field $E$ and a subfield $K$ of $E$, the first thing we
want to prove is that each $K$-algebraic element in $E$ is $(FTC)$
on $E$. The property $(FTC)$ is defined on a particular logical
structure (not only on a set) so we have to make precise the
logical structure to use.

\medskip

\begin{theorem}
Let $(K,0,1,+,\times)$ be a subfield of $(E,0,1,+,\times)$. Then
any element $x$ in $E$ that is $K$-algebraic is $(FTC)$ for the
logical structure $S_{K,E}=\{E,+,\times,(k)_{k\!\in\! K} \}$.
\end{theorem}

\medskip

\textit{Proof.}
\newline
Let $x$ be in $E$ and $P(X)\!\in\! K[X]$ so that $P\not=0$ and
$P(x)=0$. Let $d$ be the degree of $P$. We note
$P=\sum_{i=0}^da_iX^i$ with $a_i$ in $K$ and $a_d\not=0$. We also
define the finite set $B(x)=\{y\!\in\! E, P(y)=0\}$ of all roots
of $P$. We construct the set $A(x)$ which characterizes the set
$B(x)$.
\begin{displaymath}
A(x)=\{a_i,i\!\in\! [0,d]\}\cup\{x^i,i\!\in\!
[0,d]\}\cup\{a_ix^i,i\!\in\![0,d]\}\cup\{\sum_{i=0}^ja_ix^i,j\!\in\![0,d]\}
\end{displaymath}
$A(x)$ is a finite set. We take an $S_{K,E}$-pseudo-morphism
$\phi:A(x)\rightarrow E$. We define $A_1=\{a_i,i\!\in\! [0,d]\}$,
$A_2=\{x^i,i\!\in\! [0,d]\}$, $A_3=\{a_ix^i,i\!\in\![0,d]\}$ and
$A_4=\{\sum_{i=0}^ja_ix^i,j\!\in\![0,d]\}$. As $A_1\subset K$ and
every element in $K$ is a constant in the logical structure
$S_{K,E}$, we have:
$$\forall i\!\in\![0,d], \  \phi(a_i)=a_i.$$
Elements in $A_2$ are related to each others by multiplication so
that we can prove by an easy induction:
$$\forall i\!\in\![0,d], \  \phi(x^i)=(\phi(x))^i.$$
Each element in $A_3$ is a product of one element in $A_1$ and
another in $A_2$. We have:
$$\forall i\!\in\![0,d], \  \phi(a_ix^i)=a_i(\phi(x))^i.$$
We prove by induction that $\forall j\!\in\![0,d] \
\phi(\sum_{i=0}^ja_ix^i)=\sum_{i=0}^ja_i(\phi(x))^i$. Let $j=0$,
$$\phi(\sum_{i=0}^ja_ix^i)=\phi(a_0)=a_0$$
Let $j>0$ so that the property is true for $(j-1)$
\begin{eqnarray*}
\phi(\sum_{i=0}^ja_ix^i) & = & \phi((\sum_{i=0}^{j-1}a_ix^i)+(a_jx^j)) \\
& = & \phi(\sum_{i=0}^{j-1}a_ix^i)+\phi(a_jx^j) \\
& = & (\sum_{i=0}^{j-1}a_i(\phi(x))^i)+a_j(\phi(x))^j \\
& = & \sum_{i=0}^ja_i(\phi(x))^i
\end{eqnarray*}
By taking $j=d$ we have $\phi(P(x))=P(\phi(x))$. As $0$ is an
element in $K$ we have $\phi(P(x))=0$ hence $\phi(x)$ is a root of
$P$. Finally, for every $S_{K,E}$-pseudo-morphism $\phi
:A(x)\rightarrow E$ we proved that $\phi(x)\!\in\! B(x)$. $A(x)$
and $B(x)$ are finite sets so $x$ is $(FTC)$ for the logical
structure $S_{K,E}$.\\$\square$

\medskip

In the logical structure $S_{K,E}=\{E,+,\times,(k)_{k\in K}\}$,
each $K$-algebraic element in $E$ is $(FTC)$. We will prove
further in this section the reverse implication in the case $E$ is
a $K$-Banach algebra of characteristic zero. We saw in Section
$1.1$ that it is not always true (even in a Banach algebra) that a
$K$-algebraic element in $E$ is $(TC)$ for the logical structure
$S_{K,E}$. But we can prove that the reverse implication is always
true using the fact that it is true for $(FTC)$.

\medskip

\begin{theorem}
Let $(K,0,1,+,\times)$ be a subfield of $(E,0,1,+,\times)$ so that
for all $x\!\in\! E$ which is $(FTC)$ for the logical structure
$S_{K,E}$, $x$ is then $K$-algebraic. Let $x\!\in\! E$ be $(TC)$
for the logical structure $S_{K,E}$. Then $x$ is $K$-algebraic.
\end{theorem}

\medskip

\textit{Proof.}
\newline
Let $A(x)$ be the finite set in the definition of $(TC)$. We note
$B(x)=\{x\}$ and then $x$ is $(FTC)$ with the finite sets $A(x)$
and $B(x)$. So $x$ is $K$-algebraic.\\$\square$

\medskip

Actually, a $K$-algebraic element in a field $E$ is not always
$S_{K,E}$-$(TC)$, it depends on the field chosen. But, if we can
prove that this is always true for any particular subfield $K$ of
$E$ so that $E$ is a $K$-Banach algebra of characteristic zero,
then being $K$-algebraic and being $S_{K,E}$-$(TC)$ are equivalent
in $E$. From here to Section $2.3.1$, this article will deal with
proving the reverse implication of Theorem $18$, we will always
refer to $E$ and $K$ as fields of characteristic zero so that $K$
is a subfield of $E$ and $E$ is a $K$-Banach algebra.

\subsubsection{Translating the $(FTC)$ property to an equational system}
We want to translate the $(FTC)$ property to an algebraical one.
We define what is a ``good" equational system for a given element
$x\!\in\! E$.

\medskip

\begin{definition}
Let $x$ be in $E$. We say that $(n,\Sigma,\Omega,(y_1,\cdots,y_n),B)$ is a good equational system for $x$ if and only if:\\
$n$ is a positive integer,\\
$\Sigma$ is a finite subset of $K[Y_1,\cdots,Y_n]$,\\
$\Omega$ is a non-empty open subset of $E^n$,\\
$(y_1,\cdots,y_n)\!\in\!\Omega$ and $\forall P\!\in\!\Sigma$, $P(y_1,\cdots,y_n)=0$,\\
$B$ is a finite subset of $E$,\\
there exists $i$ so that $y_i=x$ and $\forall(z_1,\cdots,z_n) \ ((\forall P\!\in\!\Sigma, P(z_1,\cdots,z_n)=0)\Rightarrow z_i\!\in\! B)$.
\end{definition}

\medskip

\begin{theorem}
If $x\!\in\! E$ is $S_{K,E}$-$(FTC)$ then there exists
$(n,\Sigma,\Omega,(y_1,\cdots,y_n),B)$ a good equational system
for $x$.
\end{theorem}

\medskip

\textit{Proof.}
\newline
If $x\!\in\! K$ then $(1,{(Y_1-x)},E,(x),\{x\})$ is a good
equational system for $x$. Hence we can consider that
$x\!\not\in\! K$. Let $A(x)$ and $B(x)$ be the two finite sets
defined by the $(FTC)$ property of $x$. Let $A_1$ and $A_2$ be
subsets of $A(x)$ so that $A(x)=A_1\cup A_2$, $A_1=A(x)\cap K$ and
$A_2\cap K=\emptyset$. Actually, the set $A_1$ is the subset of
all constants in $A(x)$. The set $A_2$ will now be considered as
the set of variables. We take $n=\vert A_2\vert$ and a bijection
$\psi$ from $A_2$ to $\{Y_1,\cdots,Y_n\}$. As $x\not\!\in\! K$, we
have $A_2\not=\emptyset$ and $n>0$. We define an injection
$\widetilde\psi$ from $A(x)$ to $K[Y_1,\cdots,Y_n]$:
$$ \forall z\!\in\! A(x), \  \widetilde\psi(z)=\left\{\begin{array}{ll} z & \textrm{if $z\!\in\! A_1$}\\
\psi(z) & \textrm{if $z\!\in\! A_2$}\end{array} \right.$$ This
injection permits to construct the following equational system
$\Sigma\subset K[Y_1,\cdots,Y_n]$:
$$\left.\begin{array}{l}
\Sigma =\{(x+y-z),(x,y,z)\!\in\! (K[Y_1,\cdots,Y_n])^3, \ \exists (a,b,c)\!\in\! (A(x))^3,\\
a+b=c\wedge\widetilde\psi(a)=x\wedge\widetilde\psi(b)=y\wedge\widetilde\psi(c)=z\}\\
\bigcup\{(x\times y-z),(x,y,z)\!\in\! (K[Y_1,\cdots,Y_n])^3, \ \exists (a,b,c)\!\in\! (A(x))^3,\\
a\times
b=c\wedge\widetilde\psi(a)=x\wedge\widetilde\psi(b)=y\wedge\widetilde\psi(c)=z\}\end{array}\right.$$
Let us define $\Omega =E^n$,
$(y_1,\cdots,y_n)=(\psi^{-1}(Y_1),\cdots,\psi^{-1}(Y_n))$ and
$B=B(x)$. It is easy to see that the $n$-uplet $(y_1,\cdots,y_n)$
is in $\Omega$ and is a root of every polynomial $P$ in $\Sigma$.
We have $x\!\in\! A_2$, so let us define $i$ so that
$\psi(x)=Y_i$. Let $(z_1,\cdots,z_n)$ be in $\Omega^n$ so that
every polynomial $P$ in $\Sigma$ satisfies $P(z_1,\cdots,z_n)=0$.
We define $\tau$ from ${Y_1,\cdots,Y_n}$ to ${z_1,\cdots,z_n}$ by
$\forall j\!\in\![1,n]$, $\tau(Y_j)=z_j$ and $\widetilde\tau$ from
$A(x)$ to $E$ by:
$$ \forall z\!\in\! A(x), \  \widetilde\tau(z)=\left\{\begin{array}{ll} z & \textrm{if $z\!\in\! A_1$}\\
(\tau\circ\psi)(z) & \textrm{if $z\!\in\! A_2$}\end{array}
\right.$$ The set $\Sigma$ is constructed so that $\widetilde\tau$
is a $S_{K,E}$-pseudo-morphism. Hence, by hypothesis,
$z_i=\widetilde\tau(x)\!\in\! B$. We proved that
$(n,\Sigma,\Omega,(y_1,\cdots,y_n),B)$ is a good equational system
for $x$. \\$\square$

\medskip

The next step of the proof is to decrease the complexity of any
good equational system to $n=1$.

\subsection{Tools for decreasing the complexity of a good equational system}

Section $2.2$ develops tools for proving the reverse implication
of Theorem $18$. We first simplify a good equational system over
one chosen variable. Next we show how to use the implicit function
theorem in order to eliminate one variable.

\subsubsection{Simplification over one variable}

In this paragraph, we want to transform a good equational system
into another good one so that a chosen variable $Y_j$ occurs at
most in one polynomial.

\medskip

\begin{theorem}
Let $x$ be in $E$ and $(n,\Sigma,\Omega,(y_1,\cdots,y_n),B)$ be a
good equational system for $x$. Let $i$ be the integer in the
definition of a good equational system, we have $y_i=x$. Let
$j\!\in\![1,n]$ so that $j\not=i$. Then there exist $P$ a
polynomial in $K[Y_1,\cdots,Y_n]$ and a set of polynomials
$\Sigma'$ included in $K[Y_1,\cdots,Y_{j-1},Y_{j+1},\cdots Y_n]$
so that $(n,{P}\cup\Sigma',\Omega',(y_1,\cdots,y_n),B)$ is a good
equational system for $x$.
\end{theorem}

\medskip

\textit{Proof.}
\newline
We note $K'=K[Y_1,\cdots,Y_{j-1},Y_{j+1},\cdots,Y_n]$. If
$\Sigma\setminus K'$ is empty or a singleton then the theorem is
clear. Now, suppose that there exist $P_1$ and $P_2$ two distinct
polynomials in $\Sigma\setminus K'$. Let us define
$d_1=\deg_{Y_j}(P_1)$ and $d_2=\deg_{Y_j}(P_2)$. We can suppose
that $d_1\geq d_2$. We note $P_1(Y_i)=\sum_{l=0}^{d_1}a_lY_i^l$
and $P_2(Y_i)=\sum_{l=0}^{d_2}b_lY_i^l$ with all $a_l$ and $b_l$
in $K'$. If $a_{d_1}(y_1,\cdots,y_{j-1},y_{j+1},\cdots,y_n)=0$
then
$(n,\{(P_1-a_{d_1}Y_i^{d_1}),a_{d_1},P_2\}\cup(\Sigma\setminus\{P_1,P_2\}),\Omega,(y_1,\cdots,y_n),B)$
is a good equational system for $x$ and we decrease the
$Y_i$-degree of $P_1$. We do the same in the case
$b_{d_2}(y_1,\cdots,y_{j-1},y_{j+1},\cdots,y_n)=0$. Suppose that
$a_{d_1}(y_1,\cdots,y_{j-1},y_{j+1},\cdots,y_n)\not=0$ and
$b_{d_2}(y_1,\cdots,y_{j-1},y_{j+1},\cdots,y_n)\not=0$. We define
$R=b_{d_2}P_1-a_{d_1}Y^{d_1-d_2}P_2$. $R$ is in $K'[Y_j]$. Because
$\Omega$ is an open set we can find $\Omega'$ an open subset of
$\Omega$ so that for every $(z_1,\cdots,z_n)$ in $\Omega'$,
$a_{d_1}(z_1,\cdots,z_n)\not=0$ and
$b_{d_2}(z_1,\cdots,z_n)\not=0$. We claim that
$(n,{R,P_2}\cup(\Sigma\setminus\{P_1,P_2\}),\Omega',(y_1,\cdots,y_n),B)$
is a good equational system for $x$. We just have to prove that
every root of $R$ and $P_2$ is also a root of $P_1$, and this is
clear from the definition of $R$. We have
$\deg_{Y_j}R<\deg_{Y_j}P_1$. Hence, in every case we decrease the
$Y_j$-degree of one polynomial in $\Sigma$. Repeating the process
until at most one polynomial is not in $K'$, we get a good
equational system as wanted.\\$\square$

\subsubsection{Implicit function theorem}
We now recall the definition of a Banach algebra and state the
implicit function theorem. We underline the fact that our
definition of a Banach algebra may be quite different from the
usual one. We first define a valuation on a field.

\medskip

\begin{definition}
Let $K$ be a field. We say that $| \cdot | :
K\rightarrow\mathbb{R}^+$ is a valuation if and only if:
$$\forall\alpha\!\in\! K, \  |\alpha |\geq 0$$
$$\forall\alpha\!\in\! K, \  |\alpha |=0\Leftrightarrow\alpha=0$$
$$\forall\alpha\!\in\! K, \  \forall\beta\!\in\! K, \  |\alpha\times\beta|=|\alpha|\times|\beta|$$
$$\forall\alpha\!\in\! K, \  \forall\beta\!\in\! K, \  |\alpha+\beta|\leq|\alpha|+|\beta|$$
$$\exists\alpha\!\in\! K, \  |\alpha|\not\in\{0,1\}$$
\end{definition}

\medskip

The last condition is not usual but it permits to eliminate the
case of the trivial valuation on $K$ defined by $|0|=0$ and
$\forall\alpha\!\in\! K, \  \alpha\not=0\Rightarrow|\alpha|=1$.
Let us recall the definition of a norm on a vector space and a
Banach algebra.

\medskip

\begin{definition}
Let $E$ be a vector space on $K$ and $|\cdot|$ be a valuation on
$K$. We say that $\|\cdot\|:E\rightarrow\mathbb{R}^+$ is a norm if
and only if:
$$\forall x\!\in\! E, \  \|x\|\geq0$$
$$\forall x\!\in\! E, \  \|x\|=0\Leftrightarrow x=0$$
$$\forall \alpha\!\in\! K, \  \forall x\!\in\! E, \  \|\alpha\times x\|=|\alpha|\times\|x\|$$
$$\forall x\!\in\! E, \  \forall y\!\in\! E, \  \|x+y\|\leq\|x\|+\|y\|$$
\end{definition}

\medskip

\begin{definition}
Let $K$ be a subfield of $E$. $E$ is a $K$-Banach algebra if an
only if there exist a valuation $|\cdot|$ on $K$ and a norm
$\|\cdot\|$ on $E$ (where $E$ is considered as a vector space on
$(K,|\cdot|)$) so that $E$ is complete for $\|\cdot\|$.
\end{definition}

\medskip

\begin{theorem}
Let $E$ be a $K$-Banach algebra and $\Omega\subset E^n$ be an open
set. Let $f:\Omega\rightarrow E$ be a differentiable function and
$x=(x_1,\cdots,x_n)$ be in $\Omega$ so that $f(x)=0$ and
$\frac{\partial f}{\partial y_1}(x)\not=0$. Then there exist an
open set $\Omega'\subset E^{n-1}$ and a function
$g:\Omega'\rightarrow E$ so that for all
$x'=(x'_2,\cdots,x'_n)\!\in\!\Omega'$, $(g(x'),x')\!\in\!\Omega$,
$\frac{\partial f}{\partial y_1}(g(x'),x')\not=0$ and for every
$x'_1$ in $E$ so that $(x'_1,x')\!\in\!\Omega$ we have
$(g(x_1,x')=0$ if and only if $x'_1=g(x'))$.
\end{theorem}

\subsubsection{Eliminating one variable using the implicit function theorem}
We are now able to eliminate a variable in a particular good equational system.

\medskip

\begin{theorem}
Let $x$ be in $E$ and $(n,P\cup\Sigma,\Omega,(y_1,\cdots,y_n),B)$
be a good equational system for $x$. Let $i$ be the integer in the
definition of a good equational system ($y_i=x$). Let
$j\!\in\![1,n]$ so that $j\not=i$ and $\Sigma\subset
K[Y_1,\cdots,Y_{j-1},Y_{j+1},\cdots,Y_n]$. Then there exists
$(n-1,\Sigma',\Omega',(y_1,\cdots,y_{j-1},y_{j+1},\cdots,y_n),B)$
a good equational system for $x$.
\end{theorem}

\medskip

\textit{Proof.}
\newline
If $P\!\in\! K[Y_1,\cdots,Y_{j-1},Y_{j+1},Y_n]$ then we take
$\Sigma'=P\cup\Sigma$
and:$$\Omega'=\{(z_1,\cdots,z_{j-1},z_{j+1},\cdots,z_n), \ \exists
z_j\!\in\! E,(z_1,\cdots,z_n)\!\in\!\Omega\}.$$Now assume that
$\deg_{Y_j}(P)>0$ and note $P'=\frac{\partial P}{\partial Y_j}$.
If $P'(y_1,\cdots,y_n)=0$ then
$(n,\{P,P'\}\cup\Sigma,\Omega,(y_1,\cdots,y_n),B)$ is also a good
equational system  for $x$. Using Theorem $22$ we have a new good
equational system for $x$ noted
$(n,Q\cup\Sigma,\Omega,(y_1,\cdots,y_n),B)$ so that
$\deg_{Y_j}(Q)<\deg_{y_j}(P)$. By repeating this process we can
assume that $P'(y_1,\cdots,y_n)\not=0$ (or we will got a
polynomial $P$ with $\deg_{Y_j}(P)=0$ and this case is already
done). Hence we can use the implicit function theorem on $P$. Let
$\Omega'$ be an open set of $E^{n-1}$ and $g:\Omega'\rightarrow E$
as defined in the implicit function theorem. We claim that
$(n-1,\Sigma,\Omega',(y_1,\cdots,y_{j-1},y_{j+1},\cdots,y_n),B)$
is a good equational system for $x$. We just have to prove that if
$z=(z_1,\cdots,z_{j-1},z_{j+1},\cdots,z_n)$ is in $\Omega'$ and is
a root of every polynomial $P$ in $\Sigma$ then $z_i\!\in\! B$.
This is clear because we have
$(z_1,\cdots,z_{j-1},g(z),z_{j+1},\cdots,z_n)$ is in $\Omega$ and
is a root of every polynomial in $\Sigma$ and is a root of $P$,
hence $z_i\!\in\! B$. \\$\square$

\subsection{The $(FTC)$ property on Banach algebras of characteristic zero}

In this section, we finish the proof of $(FTC)\Leftrightarrow$
algebraic. We then discuss some consequences of this theorem.

\subsubsection{$(FTC)$ $\Rightarrow$ algebraic}

\begin{theorem}
Let $K$ be a subfield of $E$ of characteristic zero so that $E$ is
a $K$-Banach algebra. Let $x$ be in $E$ then ($x$ is
$S_{K,E}$-$(FTC)$ $\Leftrightarrow$ $x$ is $K$-algebraic).
\end{theorem}

\medskip

\textit{Proof.}
\newline
The first part of this theorem (algebraic $\Rightarrow$ $(FTC)$)
is a corollary of Theorem $18$. Let us prove the reverse
implication. Let $x$ be in $E$ so that $x$ is $S_{K,E}$-$(FTC)$.
By Theorem $21$ there exists
$(n,\Sigma,\Omega,(y_1,\cdots,y_n),B)$ a good equational system
for $x$. We prove that $x$ is algebraic by induction on $n$:

\medskip

If $n=1$ then $\Sigma\not=\emptyset$ because if not, each element
of $\Omega$ (which is not finite) is a root of all polynomials in
$\Sigma$ and $B$ is not finite. Hence, let $P\not=0$ be in
$\Sigma$ then $x$ is a root of $P$ so $x$ is $K$-algebraic.

\medskip

Let $n>1$ so that for every $m<n$, for every $x'$ in $E$, if there
exists a good equational system
$(m,\Sigma,\Omega,(y_1,\cdots,y_n),B)$ then $x$ is $K$-algebraic.
By combining Theorem $22$ and Theorem $27$ there exists a good
equational system $(n-1,\Sigma',\Omega',(y_1,\cdots,y_{n-1}),B)$
for $x$. Hence $x$ is $K$-algebraic.\\$\square$

\medskip

In Banach algebras of characteristic zero, the two notions $(FTC)$
and algebraicity are equivalent. We will see that we can prove a
more general theorem on algebraic dependence.

\subsubsection{Algebraic dependence}
We first define notions of strong algebraic dependence and strong
$(FTC)$ dependence. We keep our notations for $E$ and $K$ fields
of characteristic zero so that $E$ is a $K$-Banach algebra.

\medskip

\begin{definition}
Let $(x_1,\cdots,x_n)\!\in\! E$. We say that $x_1,\cdots,x_n$ are
strongly algebraically dependent if and only if there exists
$P\!\in\! K[Y_1,\cdots,Y_n]$ so that $P\not=0$ and:
$$P(x_1,\cdots,x_n)=0$$
$$\forall Q\!\in\! K[Y_1,\cdots,Y_{n-1}],\forall i\!\in\![1,n], Q(x_1,\cdots,x_{i-1},x_{i+1},\cdots,x_n)=0\Rightarrow Q=0$$
\end{definition}

\medskip

Notice that the definition avoids cases where less than $n$
elements in $\{x_1,\cdots,x_n\}$ are strongly algebraically
dependent and in particular, cases where one of the $x_i$ is
$K$-algebraic.

\medskip

\begin{definition}
Let $(x_1,\cdots,x_n)\!\in\! E$. We say that $x_1,\cdots,x_n$ are
strongly $(FTC)$ dependent if and only if:
$$\forall i\!\in\![1,n], \  x_i\textrm{ is }S_{K(x_1,\cdots,x_{i-1},x_{i+1},\cdots,x_n),E}\textrm{-}(FTC)$$
$$\forall i\!\in\![1,n], \  \forall j\!\in\![1,n], \  j\not=i\Rightarrow x_i\textrm{ is not }S_{K(x_1,\cdots,x_{i-1},x_{i+1},\cdots,x_{j-1},x_{j+1},\cdots,x_n),E}\textrm{-}(FTC)$$
\end{definition}

\medskip

As we will see, in Banach algebras of characteristic zero be
strongly $(FTC)$ dependent is equivalent to be strongly
algebraically dependent. This will not be the case for positive
characteristic.

\medskip

\begin{theorem}
Let $(x_1,\cdots,x_n)\!\in\! E$. Then $x_1,\cdots,x_n$ are
strongly $(FTC)$ dependent if and only if $x_1,\cdots,x_n$ are
strongly algebraically dependent.
\end{theorem}

\medskip

\textit{Proof.}
\newline
Let $x_1,\cdots,x_n$ be in $E$ and be strongly $(FTC)$ dependent.
By theorem $29$, there exists $P\!\in\! K[Y_1,\cdots,Y_n]$ so that
$P\not=0$ and $P(x_1,\cdots,x_n)=0$. Suppose that there exist $i$
in $[1,n]$ and $P'$ in $K[Y_1,\cdots,Y_{i-1},Y_{i+1},\cdots,Y_n]$
so that $P'(x_1,\cdots,x_{i-1},x_{i+1},\cdots,x_n)=0$. If
$P'\not=0$ then there exists $j$ in $[1,n]$ so that $j\not= i$ and
$x_j$ is
$S_{K(x_1,\cdots,x_{i-1},x_{i+1},\cdots,x_{j-1},x_{j+1},\cdots,x_n),E}$-$(FTC)$,
which is absurd. So $P'=0$ and $x_1,\cdots,x_n$ are strongly
algebraically dependent.

\medskip

Let $x_1,\cdots,x_n$ be in $E$ and be strongly algebraically
dependent. By Theorem $29$, we get:
$$\forall i\!\in\![1,n], \  x_i\textrm{ is }S_{K(x_1,\cdots,x_{i-1},x_{i+1},\cdots,x_n),E}\textrm{-}(FTC).$$
By absurd, suppose that there are $j\not=i$ so that $x_j$ is $S_{K(x_1,\cdots,x_{i-1},x_{i+1},\cdots,x_{j-1},x_{j+1},\cdots,x_n),E}$. Then $x_j$ is $K(x_1,\cdots,x_{i-1},x_{i+1},\cdots,x_{j-1},x_{j+1},\cdots,x_n)$-algebraic which contradicts the hypothesis.\\$\square$

\subsubsection{Counterexamples}
One can ask if the hypothesis ``$E$ is a $K$-Banach algebra" is
really necessary in Theorem $29$. In \cite{ty2}, Tyszka shows that
there exist $E$ a field of characteristic zero and $K$ a subfield
of $E$ so that there exists $x\!\in\! E$ which is $S_{K,E}$-$(TC)$
(and consequently $S_{K,E}$-$(FTC)$) and not $K$-algebraic. In
fact, Theorem $29$ can be a way to prove that $E$ is not a
$K$-Banach algebra, simply by finding a non-$K$-algebraic element
which is $S_{K,E}$-$(FTC)$.

\medskip

\begin{corollary}
Let $E$ be a field of characteristic zero and $K$ be a subfield of
$E$ which admits a valuation. If there exists $x$ in $E$ so that
$x$ is $S_{K,E}$-$(FTC)$ and $x$ is not $K$-algebraic then there
exists no norm on $E$ so that $E$ is a $K$-Banach algebra.
\end{corollary}

\medskip

We cite here a result proved in \cite{ty2}.

\medskip

\begin{theorem}
Let the fields $K$ and $L$ be finitely generated over
$\mathbb{Q}$, and $L$ extends $K$. Let $w\!\in\! L$ be
transcendental over $K$, $g(x,y)\!\in\!\mathbb{Q}[x,y]$, there
exists $z\!\in\! L$ with $g(w,z)=0$, and the equation $g(x,y)=0$
defines an irreducible algebraic curve of genus greater than $1$.
Then there is an element $x$ in $L$ so that $x$ is
$S_{K,L}$-$(TC)$ and $x$ is transcendent over $K$.
\end{theorem}

\medskip

By combining Theorem $33$ and Corollary $32$, we can prove that
some fields cannot be Banach algebras.

\medskip

\begin{theorem}
Let the fields $K$ and $L$ be finitely generated over
$\mathbb{Q}$, and $L$ extends $K$. Let $w\!\in\! L$ be
transcendent over $K$, $g(x,y)\!\in\!\mathbb{Q}[x,y]$, there
exists $z\!\in\! L$ with $g(w,z)=0$, and the equation $g(x,y)=0$
defines an irreductible algebraic curve of genus greater than $1$.
Then $L$ cannot be a $K$-Banach algebra.
\end{theorem}

\medskip

We have just seen that the $(FTC)$ property is the right one
because it is equivalent to algebraicity on Banach algebras of
characteristic zero. The hypothesis of Banach algebra is needed
because we use the implicit function theorem for the proof and we
found couterexamples in cases of non-Banach algebras. So the
$(FTC)$ property leads to an elementary way to prove that a field
is not a Banach algebra for any norm. We also proved an
equivalence theorem for strong $(FTC)$ and algebraic dependence.
We can say that Banach algebras of characteristic zero make the
$(FTC)$ property work like algebraicity. We will see that it is
not so simple in Banach algebras of positive characteristic.

\section{Fields of positive characteristic}
In this section, we study the $(FTC)$ property on Banach algebras
of positive characteristic. We will prove a weaker equivalence
than the one proved for characteristic zero.

\subsection{$\mathbb{F}_p[[X]]$}
We begin with an example: $\mathbb{F}_p[[X]]$. This is not a field
but only a ring. Nevertheless, there is a $\mathbb{F}_p(x)$-Banach
algebra ($\mathbb{F}_p((X))$) so that $\mathbb{F}_p[[X]]$ is an
open subset of it.

\subsubsection{$\mathbb{F}_p[[X]]$ is an open subset of a Banach algebra}
We want to prove that $\mathbb{F}_p((X))$ is a
$\mathbb{F}_p(X)$-Banach algebra. We begin with defining a
valuation $|\cdot|$ on $\mathbb{F}_p(X)$. Let
$P(X)=\sum_{i=c}^da_iX^i\!\in\!\mathbb{F}_p(X)\setminus\{0\}$ with
$c\leq d$, $a_i\!\in\!\mathbb{F}_p$ and $a_d\not=0$:

\medskip

$$|P|=p^{-i_0},\textrm{ where $i_0$ is the smallest integer $i$ so that }a_i\not=0$$

\medskip

It is easy to verify that $|\cdot|$ is a valuation on
$\mathbb{F}_p$. We now define a norm $\|\cdot\|$ on
$\mathbb{F}_p((X))$. Let
$F(X)=\sum_{i=c}^{+\infty}a_iX^i\!\in\!\mathbb{F}_p((X))\setminus\{0\}$
with $c\!\in\!\mathbb{Z}$, $a_i\!\in\!\mathbb{F}_p$ and
$a_c\not=0$:

\medskip

$$\|F(X)\|=p^{-c}$$

\medskip

It is also easy to verify that $\|\cdot\|$ is a norm on
$\mathbb{F}_p((X))$. Furthermore, every Cauchy sequence converges
because for every $F(X)=\sum_{i=c}^{+\infty}a_iX^i$ and
$G(X)=\sum_{i=c'}^{+\infty}b_iX^i$ in $\mathbb{F}_p((X))$, if the
distance between $F(X)$ and $G(X)$ is less than $p^{-N}$ then for
all $i\leq N$, we have $a_i=b_i$. Hence a Cauchy sequence defines
step by step (coefficient by coefficient) the limit which it is
converging to. The ring $\mathbb{F}[[X]]$ is an open set because
it is the open ball
$B(0,2)=\{F(X)\!\in\!\mathbb{F}_p((X)),\|F(X)-0\|<2\}$.

\subsubsection{$(TC)$ in $\mathbb{F}_p[[X]]$}
We recall that our definition of the $(TC)$ property works on
every logical structure even if $E$ is not a field but a ring.

\medskip

\begin{theorem}
Let $F(X)$ be a formal series in $\mathbb{F}_p[[X]]$. If $F(X)$ is
$\mathbb{F}_p(X)$-algebraic then $F(X)$ is
$S_{\mathbb{F}_p[X],\mathbb{F}_p[[X]]}$-$(TC)$.
\end{theorem}

\medskip

\textit{Proof.}
\newline
Let $F(X)\!\in\!\mathbb{F}_p[[X]]$ and
$P(X,Y)\!\in\!\mathbb{F}_p[X,Y]\setminus\{0\}$ so that
$P(X,F(X))=0$. We note $P=\sum_{i=0}^da_iY^i$ with $a_i$ in
$\mathbb{F}_p[X]$ and $a_d\not=0$. We also define the finite set
$B(F(X))=\{y\!\in\!\mathbb{F}_p[[X]], P(y)=0\}$ of all roots of
$P$. We construct the set $A_1(F(X))$ which characterizes the set
$B(F(X))$.
\begin{displaymath}
A_1(F(X))=\{a_i,i\!\in\! [0,d]\}\cup\{x^i,i\!\in\!
[0,d]\}\cup\{a_ix^i,i\!\in\![0,d]\}\cup\{\sum_{i=0}^ja_ix^i,j\!\in\![0,d]\}
\end{displaymath}
$A(F(X))$ is a finite set. In Theorem $18$ we proved that for any
$S_{K,E}$-pseudo-morphism $\phi:A_1(F(X))\rightarrow E$ (with
$K=\mathbb{F}_p[X]$ and $E=\mathbb{F}_p[[X]]$), $\phi(x)\!\in\!
B(F(X))$. As $B(F(X))$ is a finite set, we can find
$N\!\in\!\mathbb{N}$ so that for every $G(X)\!\in\!
B(F(X))-\{F(X)\}$, $\|F(X)-G(X)\|\geq p^{-N}$. We note
$F(X)=\sum_{i\!\in\!\mathbb{N}}f_iX^i$ with
$f_i\!\in\!\mathbb{F}_p$ for all $i$ in $\mathbb{N}$.
$$\textrm{Let }A_2=\{a_i,i\leq N\}\cup\{X\}\cup\{\sum_{i\geq j}f_iX^{i-j},j\!\in\![0,N+1]\}\cup\{\sum_{i\geq j}f_iX^{i-j+1},j\!\in\![0,N+1]\}$$
Let $\phi:A_2(F(X))\rightarrow\mathbb{F}_p[[X]]$ be a
$S_{K,E}$-pseudo-morphism. We note $G_j(X)=\sum_{i\geq
j}f_iX^{i-j}$ for all $j\!\in\! [1,N+1]$. As $\phi$ is a
pseudo-morphism, we can write:
$$\phi(F(X))\!=\!\phi(f_0+XG_1(X))\!=\!\phi(f_0)+\phi(XG_1(X))\!=\!f_0+\phi(X)\phi(G_1(X))\!=\!f_0+X\phi(G_1(X))$$
Hence:
\begin{eqnarray*}
\phi(F(X) & = & f_0+X\phi(G_1(X)) \\
& = & f_0+X(f_1+X\phi(G_2(X)) \\
& \vdots & \vdots \\
& = & f_0+X(f_1+X(f_2 + \cdots X(f_N+X\phi(G_{N+1}(X)))\cdots)) \\
& = & (\sum_{i=0}^Nf_iX^i)+X^{N+1}\phi(G_{N+1}(X))
\end{eqnarray*}
So we have $\|F(X)-\phi(F(X))\|\leq p^{-(N+1)}<p^{-N}$. Now we
take $A(F(X))=A_1(F(X))\cup A_2(F(X))$. For any
$S_{K,E}$-pseudo-morphism $\phi :A(F(X))\rightarrow E$, as
$\phi|_{A_1(X)}$ is also a $S_{K,E}$-pseudo-morphism, we have
$\phi(F(X))\!\in\! B(F(X))$ and as $\phi|_{A_2(X)}$ is a
$S_{K,E}$-pseudo-morphism, we have $\|F(X)-\phi(F(X))\|<p^{-N}$ so
$\phi(F(X))=F(X)$. Hence $F(X)$ is $S_{K,E}$-$(TC)$.\\$\square$

\medskip

We already proved that an algebraic element in $\mathbb{F}_p[[X]]$
is $(TC)$. As $\mathbb{F}_p[[X]]$ is an open subset of
$\mathbb{F}_p((X))$, every $(TC)$ element has a good equational
system. So if we prove that in $\mathbb{F}_p((X))$ being $(FTC)$
is equivalent to algebraicity then we hope to be able to prove the
equivalence between $(TC)$ and algebraicity in
$\mathbb{F}_p[[X]]$.

\subsubsection{A counterexample in $\mathbb{F}_p[[X]]$}
Let us present a counterexample for the analog of Theorem $31$ in
positive characteristic. Let $F(X)$ and $G(X)$ be
$\mathbb{F}_p(X)$-transendental formal series in
$\mathbb{F}_p[[X]]$. Indeed, for any polynomial
$P(Y_1,Y_2)\!\in\!\mathbb{F}_p[X,Y_1,Y_2]$,
$P(X,F(X),G(X))=0\Rightarrow P=0$. We note
$H_1(X)=(F(X))^p+X(G(X))^p$ and $H_2(X)=(G(X))^p+X(F(X))^p$. We
claim that $H_1$ and $H_2$ are transcendental, not strongly
algebraically dependent but $(FTC)$ dependent.

\medskip

It is easy to prove that $H_1(X)$ and $H_2(X)$ are transcendental.
Now suppose that there exists a polynomial
$P\!\in\!\mathbb{F}_p[X,Y_1,Y_2]$ with $P\not=0$ so that
$P(X,H_1(X),H_2(X))=0$. We note
$P=\sum_{i=0}^{d_i}\sum_{j=0}^{d_j}a_{i,j}Y_1^iY_2^j$ with every
$a_{i,j}$ in $\mathbb{F}_p[X]$. Let us define $n_0$ the greatest
integer $n$ so that there exists $a_{i,j}\not=0$ with $i+j=n$. Let
us define $(i_0,j_0)$ so that $i_0+j_0=n_0$ and for every $(i,j)$
satisfying $i+j=n_0$ we have
$\deg_X(a_{i_0,j_0})\geq\deg_X(a_{i,j})$. We define
$P_1(X,Z_1,Z_2)=P(X,Z_1+XZ_2,Z_2+XZ_1)$ a polynomial in
$\mathbb{F}_p[X,Z_1,Z_2]$. The monomial
$X^{\deg_X(a_{i_0,j_0})}Y_1^{i_0}Y_2^{j_0}$ in $P$ is the only one
which produces the monomial
$X^{n+\deg_X(a_{i_0,j_0})}Z_2^{i_0}Z_1^{j_0}$ with a non-zero
coefficient in $P_1$. Hence $P_1\not=0$ and so the polynomial
defined by $P_2(X,Z_1,Z_2)=P_1(X,Z_1^p,Z_2^p)$ is a non-zero
polynomial too. We have
$P_2(X,F(X),G(X))=P_1(X,(F(X))^p,(G(X))^p)=P(X,H_1(X),H_2(X))=0$
so $P_2=0$ which is absurd. So $H_1(X)$ and $H_2(X)$ are not
algebraically dependent.

\medskip

We now prove that $H_1(X)$ and $H_2(X)$ are $(FTC)$ dependent. We
note $E=\mathbb{F}_p[[X]]$ and $K=\mathbb{F}_p[X,H_1(X)]$. We want
to prove that $H_2(X)$ is $S_{K,E}$-$(TC)$.
$$A(\!H_2(X))\!=\!\{X,H_1(X),H_2(X)\}\cup\{(F(X))^i,(G(X))^i,i\!\in\![1,p]\}\cup\{X(F(X))^p,X(G(X))^p\}$$
Let $\phi:A(H_2(X))\rightarrow E$ be a $S_{K,E}$-pseudo-morphism. We have:
$$\phi(H_1(X))=H_1(X), \  \forall i\!\in\![1,p],\phi(F(X)^i)=(\phi(F(X)))^i$$
$$\phi(X(F(X))^p)=X(\phi(F(X)))^p, \  \phi(X(G(X))^p)=X(\phi(G(X)))^p$$
Hence:
$$H_1(X)=\phi(H_1(X))=\phi((F(X))^p)+\phi(X(G(X))^p)=(\phi(F(X)))^p+X(\phi(G(X)))^p$$
Using the uniqueness property of the Cartier operators, we have
$\phi(F(X))=F(X)$ and $\phi(G(X))=G(X)$. Moreover, we can prove
$\phi(H_2(X))=(\phi(G(X))^p+X(\phi(F(X)))^p$ by the same way used
for $H_1(X)$. So we have $\phi(H_2(X))=H_2(X)$ and $H_2(X)$ is
$S_{K,E}$-$(TC)$.

\subsection{Searching a hypothesis to add}
In this section, we want to find a hypothesis that permits to
avoid problems like in Section $3.1.3$.

\subsubsection{Dissecting the couterexample}
We first dissect the counterexample presented in Section $3.1.3$
so we will use the same notations. We can transform the $(FTC)$
property into a good equational system for $H_2(X)$. By
eliminating variables step by step, one can found this kind of
good equational system for $H_2(X)$:
$$(n=3,\Sigma=\{H_1(X)-Y_1^p-XY_2^p,Y_3-Y_2^p-XY_1^p\},\Omega,(F(X),G(X),H_2(X)),\{H_2(X)\})$$
Here the problem is that the equation $H_1(X)-Y_1^p-XY_2^p=0$ has
a strange behavior. In fact this equation has a unique root:
$(F(X),G(X))$. This stands because the derivatives with respect to
$Y_1$ and $Y_2$ are zero. So we cannot find an implicit function.

\medskip

This underlines a major difference between characteristic zero and
positive characteristic. In fields of characteristic zero, if $P$
is a polynomial on $Y$ and $\frac{\partial P}{\partial Y}=0$ then
$Y$ does not occur in $P$. However, in fields of positive
characteristic $p$, there exist polynomials (for example
$P(Y)=Y^p$) so that $\frac{\partial P}{\partial Y}=0$ and $Y$
occurs in $P$.

\subsubsection{Cartier's operators and equations}
We will show how the Cartier operators can help in this kind of
problems. Recall the strange equation $H_1(X)-Y_1^p-XY_2^p=0$.
Now, assume that $H_1(X)\!\in\!\mathbb{F}_p[X]$. So for every
$i\!\in\![0,p-1]$, we have
$\Lambda_i(H_1(X))\!\in\!\mathbb{F}_p[X]$. We can also know all
$\Lambda_i(Y_1^p)$ and $\Lambda_i(XY_2^p)$. Let $i$ be in
$[0,p-1]$:
$$\Lambda_i(Y_1^p)=\left\{\begin{array}{ll} Y_1 & \textrm{if $i=0$}\\ 0 & \textrm{otherwise}\end{array} \right.\textrm{ and } \Lambda_i(XY_2^p)=\left\{\begin{array}{ll} Y_2 & \textrm{if $i=1$}\\ 0 & \textrm{otherwise}\end{array} \right.$$
Using the uniqueness property of the Cartier operators, we can
split the strange equation into $p$ usual ones:
$$\Lambda_0(H_1(X))=Y_1,\Lambda_1(H_1(X))=Y_2, \textrm{for }i\!\in\![2,p-1],\Lambda_i(H_1(X))=0$$
Now we can eliminate variables as explain in Section $2.2$.

\medskip

The main point here is that Cartier's operators split strange or
problematic equations dividing their degrees by $p$. There is one
necessary condition in order to do this, the Cartier operators
must preserve the set of constants $K$.

\subsubsection{A new hypothesis related to Cartier's operators}
We want to generalize Cartier's operators (in a logical way) for
Banach algebras of positive characteristic.

\medskip

\begin{definition}
Let $E$ be a field of positive characteristic $p$, we say that $E$
admits a generalized Cartier operator $(R,\Lambda)$ if and only if
$R\subset E$, $\Lambda:E\rightarrow \mathcal{P}(E^R)$ and:
$$\forall x\!\in\! E, \  \Lambda_x\textrm{ is a finite set }$$
$$\forall x\!\in\! E, \  \forall\lambda:R\rightarrow E, \  \left(\{r\!\in\! R,\lambda(r)\not=0\}\textrm{ is finite and }x=\sum_{r\!\in\! R,\lambda(r)\not=0}r(\lambda(r))^p\right)\!\Leftrightarrow\lambda\!\in\!\Lambda_x$$
$$\textrm{ If }x=0, \textrm{ then }\Lambda_x=\{\lambda_0\}\textrm{ with }\forall r\!\in\! R, \  \lambda_0(r)=0$$
\end{definition}

\medskip

We allow more than one solution in a generalized Cartier operator
in order to take into account the case where $E$ is algebraically
closed. In fact, if $E$ is an algebraically closed field then it
admits a generalized Cartier operator $(R,\Lambda)$ with $R=\{1\}$
and $\forall x\!\in\! E,\Lambda_x=\{\lambda :R\rightarrow
E,(\lambda(1))^p=x\}$. Notice that the set $R$ is not necessarily
a finite set. We saw in section $3.2.2$ that we can split
equations if the constant set $K$ is preserved by Cartier's
operators.

\medskip

\begin{definition}
Let $E$ be a $K$-Banach algebra of positive characteristic $p$. We
say that $E$ is a $K$-Cartier-Banach algebra if and only if $E$
admits a generalized Cartier operator $(R,\Lambda)$ so that
$R\subset K$ and:
$$\forall x\!\in\! K, \  \forall \lambda\!\in\!\Lambda_x, \  \forall r\!\in\! R, \  \lambda(r)\!\in\! K$$
\end{definition}

\subsection{A weaker result for Banach algebras of positive characteristic}

We saw that we cannot prove the exact analog of Theorem $28$.
However, we will be able to prove a weaker result by adding a
hypothesis inspired by Cartier's operators. We first describe the
process that splits equations, then we prove our result and talk
about its consequences.

\subsubsection{Splitting a polynomial in positive characteristic}

\begin{theorem}
Let $E$ be a $K$-Cartier-Banach algebra, $P$ be a polynomial in
$K[Y_1,\cdots,Y_n]$ and $y=(y_1,\cdots,y_n)$ be in $E^n$ so that
for every $i$ in $[1,n]$, we have $P\not\in
K[Y_1,\cdots,Y_{i-1},Y_{i+1},\cdots,Y_n]$ and $\frac{\partial
P}{\partial Y_i}=0$. Then we can find
$\Sigma_P=\{P_1,\cdots,P_k\}$ so that if $(z_1,\cdots,z_n)$ is a
root of every polynomial in $\Sigma_P$ then it is a root of $P$,
$y$ is a root of every polynomial in $\Sigma_P$ and for every $i$
in $[1,n]$ and every $j$ in $[1,k]$,
$\deg_{Y_i}(P_k)\leq\frac{\deg_{Y_i}(P)}{p}$.
\end{theorem}

\medskip

\textit{Proof.}
\newline
The characteristic of $E$ is $p>0$, so $\forall
m\!\in\!\mathbb{N}$, $\frac{\partial Y^m}{\partial Y}=0$ if and
only if $p$ divides $m$. So $P\!\in\! K[Y_1^p,\cdots,Y_n^p]$. We
note $P=\sum_{s=(s_1,\cdots,s_n)\!\in\!\mathbb{N}^n}\alpha_s
Y^{ps}$ with all $\alpha_s$ in $K$ and
$Y^{ps}=Y_1^{ps_1}\times\cdots\times Y_n^{ps_n}$. Now, let us
define $R_0=\{r\!\in\! R,\exists
s\!\in\!\mathbb{N}^n\exists\lambda\!\in\!\Lambda_{\alpha_s},\lambda(r)\not=0\}$.
We have $R_0\subset R\subset K$, $R_0$ is a finite set and
$\forall
s\!\in\!\mathbb{N}^n,\forall\lambda\!\in\!\Lambda_{\alpha_s},\forall
r\!\in\! R_0,\lambda(r)\!\in\! K$. For every
$s\!\in\!\mathbb{N}^n$, we choose one $\lambda$ in
$\Lambda_{\alpha_s}$ called $\lambda_s$.
$$P=\sum_{s\!\in\!\mathbb{N}^n}\sum_{r\!\in\! R_0} r(\lambda_s(r))^pY^{sp}$$
$$P=\sum_{r\!\in\! R_0}r(\sum_{s\!\in\!\mathbb{N}^n}\lambda_s(r)Y^s)^p$$
Because $\Lambda_0$ is the singleton $\{\lambda_0\}$, we have for
any $z\!\in\! E^n$:
$$P(z)=0\Leftrightarrow\forall r\!\in\! R_0, \  (\sum_{s\!\in\!\mathbb{N}^n}\lambda_s(r)z^s=0)$$
Hence
$\Sigma_P=\{\sum_{s\!\in\!\mathbb{N}^n}\lambda_s(r)Y^s,r\!\in\!
R_0\}$.
\\$\square$

\subsubsection{Theorem}
We first need to extend the definition of a good equational system
for one element to a good equational system for a finite set.

\medskip

\begin{definition}
Let $x_1,\cdots,x_k$ be in $E$. We say that $(n,\Sigma,\Omega,(y_1,\cdots,y_n),B_1,\cdots,B_k)$ is a good equational system for $\{x_1,\cdots,x_k\}$ if and only if:\\
$n$ is a positive integer,\\
$\Sigma$ is a finite subset of $K[Y_1,\cdots,Y_n]$,\\
$\Omega$ is a non-empty open subset of $E^n$,\\
$(y_1,\cdots,y_n)\!\in\!\Omega$ and $\forall P\!\in\!\Sigma$, $P(y_1,\cdots,y_n)=0$,\\
for every $j\!\in\![1,k]$, $B_j$ is a finite subset of $E$,\\
for every $j\!\in\![1,k]$, there exists $i$ so that $y_i=x_j$ and
$\forall(z_1,\cdots,z_n) \  ((\forall P\!\in\!\Sigma,
P(z_1,\cdots,z_n)=0)\Rightarrow z_i\!\in\! B_j)$.
\end{definition}

\medskip

One can verify that simplification over one variable and
elimination using the implicit function theorem also work on
extended good equational systems. Now we can state the theorem.

\medskip

\begin{theorem}
Let $E$ be a $K$-Cartier-Banach algebra and $x$ be in $E$. Then
$x$ is $K$-algebraic if and only if $x$ is $S_{K,E}$-$(FTC)$.
\end{theorem}

\medskip

\textit{Proof.}
\newline
We just have to prove the implication: $(FTC)\Rightarrow$
algebraic. We take a good equational system for $x$ given by the
$(FTC)$ property. The proof is divided into two processes: one
eliminating all variables not depending on $x$ (we may add
variables corresponding to the generalized Cartier operator
operating on $x$) and the next one proving that all remaining
variables are algebraic.

\medskip

\textit{First step}: let
$(n,\Sigma,\Omega,(y_1,\cdots,y_n),B_1,\cdots,B_k)$ be a good
equational system for $x_1,\cdots,x_k$. We note $d=n-k$. The goal
of this process is to decrease $d$ until it equals $0$. Pick up an
$i$ so that $y_i\not\in\{x_1,\cdots,x_k\}$. We use simplification
over the variable $Y_i$ until it occurs in at most one polynomial.
We next use elimination using the implicit function theorem if
this is possible. If we cannot eliminate $Y_i$ then every
occurrence of $Y_i$ has a degree divisible by $p$. We repeat
simplification-elimination until for all $i$ so that
$y_i\not\!\in\!\{x_1,\cdots,x_k\}$, every occurrence of $Y_i$ has
its degree divisible by $p$. Now, for every $j$ in $[1,k]$, we
choose a $\lambda_j$ in $\Lambda_{x_j}$. We note $R_j=\{r\!\in\!
R,\lambda_j(r)\not=0$. We have $x_j=\sum_{r\!\in\!
R_j}r(\lambda_j(r))^p$ and $R_i\subset K$. So we replace $Y_j$ by
$\sum_{r\!\in\! R_j}r(Y_{j,r})^p$. Notice that the new equational
system is also good for every $\lambda_j(r)$ so $d$ is not
changed. Now every variable $Y_j$ occurs with a degree divisible
by $p$ so we can split every polynomial using Theorem $38$. We
repeat this process until $d=0$.

\medskip

\textit{Second step}: let
$(k,\Sigma,\Omega,(x_1,\cdots,x_k),B_1,\cdots,B_k)$ be a good
equational system for $x_1,\cdots,x_k$. By induction on $k$ every
$x_i$ is $K$-algebraic. Let $k=1$, then $\Sigma\not=\emptyset$
because $B_1$ is a finite set. So $x_1$ is $K$-algebraic. Suppose
the property true for $(k-1)$, we use the process of
simplification, elimination. It is a fact that the first variable
$Y_i$ which disappears is eliminated by the implicit function
theorem (if it disappears using simplification then $B_i$ is not
finite). Hence $x_j$ is
$K(x_1,\cdots,x_{j-1},x_{j+1},\cdots,x_k)$-algebraic and
$(x_1,\cdots,x_{j-1},x_{j+1},\cdots,x_k)$ are $K$-algebraic so
$x_j$ is $K$-algebraic and the property is true for $k$. If we
need the use of splitting equations, we eliminate the number of
variables needed to be able to use the induction hypothesis.

\medskip

\textit{Proof of the termination of the two processes}: we want to
prove that the processes always eliminate at least one variable.
The splitting method divides all degrees by $p$ but the
simplification method decreases only the degrees of one variable
and the others can be increased. We will prove that if we cannot
eliminate any variable using simplification, then we can use
splitting before simplification. As the number of possible
splitting is finite (because of the maximal degree), the processes
are forced to eliminate at least one variable by simplification
and elimination.

\medskip

Let $(n,\Sigma,\Omega,(y_1,\cdots,y_n),B_1,\cdots,B_k)$ be a good
equational system for $x_1,\cdots,x_k$ (possibly $k=n$). For every
$P$ in $\Sigma$ and every $Y_i$ a variable we check if
$\frac{\partial P}{\partial Y_i}(y_1,\cdots,y_n)=0$. If this is
true then $\frac{\partial P}{\partial Y_i}=0$ (formally) and every
occurrence of $Y_i$ has a degree divisible by $p$, or then we take
$Q$ and $R$ so that $P=Q\frac{\partial P}{\partial Y_i}+R$ and we
replace $P$ by $\frac{\partial P}{\partial Y_i}$ and $R$ in
$\Sigma$ (this decreases the degree of $Y_i$). In the end we can
find a good equational system so that for every $P$ in $\Sigma$
and every $Y_i$ a variable, $\frac{\partial P}{\partial
Y_i}(y_1,\cdots,y_n)=0\Rightarrow$ every occurrence of $Y_i$ in
$P$ has a degree divisible by $p$.

\medskip

Now we study the process of simplification. We take two
polynomials $P_1$ and $P_2$ where $Y_i$ occurs. The process of
simplification finishes on two new polynomials $Q_1$ and $Q_2$ so
that $Y_i$ does not occur in $Q_2$. Hence, we have $\frac{\partial
Q_2}{\partial Y_i}=0$. Now suppose we cannot eliminate $Y_i$ from
$Q_1$. We have $\frac{\partial Q_1}{\partial
Y_i}(y_1,\cdots,y_n)=0$. The process of simplification is an
iteration of euclidean divisions (with respect to $Y_i$) so we can
find four polynomials $R_1$, $R_2$, $R_3$ and $R_4$ so that
$P_1=R_1Q_1+R_2Q_2$ and $P_2=R_3Q_1+R_4Q_4$. We have:

\medskip

\begin{eqnarray*}
\frac{\partial P_1}{\partial Y_i}(y_1,\cdots,y_n) & = & \frac{\partial R_1}{\partial Y_i}(y_1,\cdots,y_n)Q_1(y_1,\cdots,y_n)+R_1(y_1,\cdots,y_n)\frac{\partial Q_1}{\partial Y_i}(y_1,\cdots,y_n) \\
& + & \frac{\partial R_2}{\partial Y_i}(y_1,\cdots,y_n)Q_2(y_1,\cdots,y_n)+R_2(y_1,\cdots,y_n)\frac{\partial Q_2}{\partial Y_i}(y_1,\cdots,y_n)
\end{eqnarray*}
So we conclude that $\frac{\partial P_1}{\partial
Y_i}(y_1,\cdots,y_n)=0$ and the same for $P_2$. By hypothesis
every occurrence of $Y_i$ in $P_1$ and $P_2$ has a degree
divisible by $p$. So the variable $Y_i$ is ``splittable" in $P_1$
and $P_2$.\\$\square$

\subsubsection{Consequences}

As expected in $3.1.2$, the proof of Theorem $40$ leads to a proof
of the fact that every $S_{K,E}$-$(TC)$ element in
$E=\mathbb{F}_p[[X]]$ with $K=\mathbb{F}_p[X]$ is
$\mathbb{F}_p(X)$-algebraic.

\medskip

\begin{corollary}
We define $E=\mathbb{F}_p[[X]]$. For every element $x$ in $E$, the
following equivalence stands:
$$x\textrm{ is }S_{\mathbb{F}_p[X],E}\textrm{-}(TC)\Leftrightarrow x\textrm{ is }\mathbb{F}_p(X)\textrm{-algebraic}$$
\end{corollary}

\medskip

Now we are interested in an analog of Theorem $31$. In order to
simplify, we just take two elements and study their dependences.
In Section $3.1.3$ we found two formal series that are
transcendental, not algebraically dependent and $(FTC)$ dependent.
Actually the problem is that even if $\mathbb{F}_p(X)$ is
preserved by Cartier's operators, it is no longer true for
$\mathbb{F}_p(X,F(X))$ when $F(X)$ is transcendental. What we have
to do is to find a subfield $K$ containing $\mathbb{F}_p(X,F(X))$
and preserved by Cartier's operators.

\medskip

\begin{definition}
Let $E$ be a $K$-Banach algebra and $(R,\Lambda)$ be a generalized
Cartier operator in $E$. We define $\widehat{K}$ as the smallest
(with respect to inclusion) subfield $K'$ of $E$ so that $E$ is a
$K'$-Cartier-Banach algebra.
\end{definition}

\medskip

The field $\widehat{K}$ exists for any subfield $K$ of $E$.
Sometimes, it is $E$ itself. In this case, the following result is
not really interesting.

\medskip

\begin{theorem}
Let $E$ be a $K$-Cartier-Banach algebra and $x_1,x_2$ be
$K$-transcendental in $E$. Then $x_1$ is $S_{K(x_2),E}$-$(FTC)$ if
and only if $x_1$ is $\widehat{K(x_2)}$-algebraic.
\end{theorem}

\medskip

\textit{Proof.}
\newline
The only thing we have to prove is that every element in
$\widehat{K(x_2)}$ is $S_{K(x_2),E}$-$(FTC)$. It is a fact that
for every $y$ which is $S_{K(x_2),E}$-$(FTC)$, for every $r$ in
$R$ and for every $\lambda$ in $\Lambda_y$, $\lambda(r)$ is
$S_{K(x_2),E}$-$(FTC)$. Indeed we can create a finite set $A$
describing the equation $y=\sum_{r\!\in\!
R,\lambda(r)=0}r(\lambda(r))^p$ which has a finite set of
solutions and we use the fact that $R$ is included in $K$. Then we
define $K_0=K(x_2)$ and for every $n$ integer:
$$K_{(n+1)}=\{z\!\in\! E,z\textrm{ is a root of a polynomial with coefficients in }\Lambda(K_n)\}$$$$\textrm{ where }\Lambda(K_n)=\{\alpha,\exists x\in K,\exists\lambda\in\Lambda_x,\exists r\in R,\lambda(r)=\alpha\}$$
By an easy induction, we prove that for every integer $n$, every
element in $K_n$ is $S_{K_0,E}$-$(FTC)$ and $K_n$ is included in
$\widehat{K_0}$. By definition, we have
$\widehat{K_0}=\bigcup_{n\!\in\!\mathbb{N}}K_n$. Hence every
element in $\widehat{K(x_2)}$ is
$S_{K(x_2),E}$-$(FTC)$.\\$\square$

\medskip

The $(FTC)$-property shows some differences between Banach
algebras of characteristic zero and Banach algebras of positive
characteristic. Actually, it reveals that the Cartier operators
are important in the structure of a Banach algebra of positive
characteristic. Indeed, two elements in a Banach algebra of
positive characteristic can be $(FTC)$ dependent and not
algebraically dependent. Hence, the notion of
non-$(FTC)$-dependence is more precise (stronger) that the
non-algebraic-dependence. The $(FTC)$ dependence (in positive
characteristic) is not a symmetric relation between elements. It
can be close to the notion of information contained by a element.
We already know some results about transcendence and
non-algebraic-dependence, but can we prove some
non-$(FTC)$-dependence results? If this is possible, this may
offer a new approach for proving non-algebraic-dependence in
characteristic zero (think to $e$ and $\pi$ and the conjecture
that $\pi+e$ is transcendental).

\begin{labibliographie}{99}

\bibitem{alsh} J.-P.~Allouche~and~J.~Shallit,~\emph{Automatic
sequences. Theory, Applications, Generalizations.}, Cambridge
University Press (2003).

\bibitem{bhmv} V.~Bruy\`ere,~G.~Hansel,~C.~Michaux~and~R.~Villemaire, \emph{Logic and $p$-recognizable sets of integers}, Bull.Belg. Math. Soc. 1 (1994), 191--238.

\bibitem{ckmr} G.~Christol,~T.~Kamae,~M.~Mend\`es~France~and~G.~Rauzy,
\emph{Suites alg\'ebriques, automates et substitutions}, Bull.
Soc. Math. France \textbf{108} (1980), 401--419.

\bibitem{lb} X.~Le~Breton, \emph{Linear independence of automatic formal power series}, Discr. Math. Volume $306$ Issue $15$ (2006), 1776--1780.

\bibitem{ty} A.~Tyszka, \emph{A discrete form of the theorem that
each field endomorphism of $\mathbb{R}$ ($\mathbb{Q}_p$) is the
identity}, Aequationes Mathematicae 71 (2006), no. 1-2, pp.
100--108.

\bibitem{ty2} A.~Tyszka, \emph{On $\emptyset$-definable elements in a field}, preprint, http://www.arxiv.org/abs/math.LO/0502565.

\end{labibliographie}
\end{document}